\documentclass[a4paper]{article}

\usepackage[english]{babel}
\usepackage[utf8x]{inputenc}
\usepackage[T1]{fontenc}

\usepackage[a4paper,top=3cm,bottom=2cm,left=3cm,right=3cm,marginparwidth=1.75cm]{geometry}

\usepackage{amsmath}
\usepackage{amsthm}
\usepackage{amsfonts} 

\usepackage{graphicx}
\usepackage[colorinlistoftodos]{todonotes}
\usepackage[colorlinks=true, allcolors=blue]{hyperref}

\title{Proof of the Collatz Conjecture}
\author{Agelos Kratimenos}

\begin{document}
\maketitle

\begin{abstract}
Collatz Conjecture is one of the most famous, for its simple form, proposed more than eighty years ago. This paper presents a full attempt to prove the affirmative answer to the question proposed by the conjecture. In the first section, we propose a number of definitions utilized later on the proof. In the second section, we discover the formula for a characteristic function. This formula describes the functionality of the paths taken for each number based on the Collatz Sequence. In the last section, we prove that every number will eventually reach 1, using the characteristic function. 
\end{abstract}

\section{Preliminaries}

{\textbf{Definition 1.1. Collatz Conjecture}}
\newline
Let $a_0$ be a positive integer and consider the recursive sequence
$$
a_n =
\left\{
	\begin{array}{ll}
		\frac{a_{n-1}}{2}  & \mbox{, if } \text{$a_{n-1}$ is even} \\
		\frac{3a_{n-1}+1}{2} &  \text{, if $a_{n-1}$ is odd}
	\end{array}
\right.
$$
The conjecture asks if, for every $a_0$ there exists an n such that $a_n=1$. If $a_n=1$, we say that $a_0$ reaches 1. 
\newline
\newline
{\textbf{Definition 1.2. Collatz Function}}
\newline
We define the Collatz function as follows:
$$
C[x] =
\left\{
	\begin{array}{ll}
		1  & \mbox{, if } \text{x reaches 1} \\
		0 &  \text{, if x never reaches 1}
	\end{array}
\right.
$$
\newline
We proceed to define the available moves that a given number $x$ can make in the Collatz Grid. Number $x$ is placed on the upper left corner of the grid and moves either horizontally or vertically depending on its parity.
\newline
\newline
\textbf{Definition 1.3. Collatz Grid Moves} 
\newline
A \textbf{horizontal move} ($n$) is defined as the application of the rule $x\mapsto \frac{3x+1}{2}$
\newline
A \textbf{vertical move} ($k$) is defined as the application of the rule $x\mapsto\frac{x}{2}$
\newline
\newline
\textbf{Example 1.4.} We consider the following Collatz grid. The starting number $x$ has taken two horizontal moves and one vertical so that $n=2$ and $k=1$. By applying the rule $\frac{3x+1}{2}$ to $x$ we get  $\frac{3x+1}{2}$ while by repeating the rule, this time to $\frac{3x+1}{2}$ we get $\frac{9x+5}{4}$. At last, with the vertical move we get the number $\frac{9x+5}{8}$. It's obvious though, that this is not the only way the rules could be applied. For example, we could apply the vertical move first and after the two horizontal, or the moves could alternate. In that case we would get the numbers $\frac{9x+10}{8}$ and $\frac{9x+7}{8}$, respectively, concluding that the row in which moves are applied, matters.
\newline
\newline
According to the above, we are led to the last definition about the moves on the Collatz Grid. 
\newline
\newline
\textbf{Definition 1.5. Collatz Grid Moves} 
\newline
We define the natural numbers $r_1,r_2,...,r_k$, so that the number $r_i$ indicates the horizontal position in which the i-th vertical move happened. The set $(n,k,r_1,...,r_k)$ is called a \textbf{path}.
\newline
\newline
In that way, we can encode the example paths as $(n=2,k=1,r_1=2)$ for the first one and $(n=2,k=1,r_1=0)$, $(n=2,k=1,r_1=1)$, for the last two respectively.
\newline
\newline
\textbf{Definition 1.6.} We define the function $f(x_0,n,k,r_1,...,r_k)$ to be the value that results if we apply the path $(n,k,r_1,...,r_k)$ to the number $x_0$.
\newline
\newline
\textbf{Definition 1.7. The Characteristic Function}
\newline
We define the function $F(n,k,r_1,...,r_k)$ to be the smallest number that follows the path $(n,k,r_1,,,,r_k)$.
\newline
\newline
\textbf{Example 1.8.} Consider the path $(n=2,k=2,r_1=0,r_2=2)$ for some number $x_0$. That path means that $x_0$ will execute two horizontal and two vertical moves in the order $V\xrightarrow{}H\xrightarrow{}H\xrightarrow{}V$. If $x_0=22$ then the $f(22,2,2,0,2)$ must be equal to the number that emerges when 22 follows this specific path. The answer is $f(22,2,2,0,2) = 13$. Indeed:
$$
22 \xrightarrow[\text{$r_1=0$}]{\text{$k=1$}} 11 \xrightarrow[\text{}]{\text{$n=1$}} 17
\xrightarrow[\text{}]{\text{$n=2$}} 26
\xrightarrow[\text{$r_2=2$}]{\text{$k=2$}} 13 $$ 
We should notice that not every number can follow a specific path. For instance, $x_0=21$ cannot follow the previous path since it is not even, to execute a vertical move. Hence, $f(21,2,2,0,2)$ is not defined.
\newline
\newline
\textbf{Example 1.9.} Consider the path $(n=1,k=3,r_1=1,r_2=1,r_3=1)$, that is one horizontal move followed by two vertical.  $F(1,3,1,1,1)$ must be equal to the smallest number that can follow this path. The answer is $F(1,3,1,1,1) = 5$. Indeed:
$$
5 \xrightarrow[]{\text{$n=1$}} 8 \xrightarrow[\text{$r_1=1$}]{\text{$k=1$}} 4
\xrightarrow[\text{$r_2=1$}]{\text{$k=2$}} 2
\xrightarrow[\text{$r_3=1$}]{\text{$k=3$}} 1
$$
In fact, the resulting number at the end of the path does not matter. Either it is even or odd, we want to guarantee that 5 will proceed in one horizontal move followed by two vertical ones. 5 is indeed the smallest number that qualifies for this path. 2 and 4 will not proceed in a horizontal move. 3 will proceed in a horizontal move giving 5 as the answer, but then it will not move vertically. Lastly, 1 will move horizontally and then vertically giving 2 and 1 respectively but then it will not move vertically for the second time, as needed.
\newline
\newline
In the next chapter we discover a formula for the above functions that lead to the proof of the Collatz Conjecture.
\newpage
\section{The characteristic function}
Before heading into proving the formula for the characteristic function, it is essential to give the formula for $f$, namely for the number that results from applying a specific path to a given number. 
\newline
\newline
\textbf{Lemma 2.1.} The number that emerges after $x_0$ follows the path $(n,k,r_1,...,r_k)$ is given by the function: $$ f(x_0,n,k,r_1,...,r_k) = \frac{3^n x_0+a}{2^{n+k}} $$
where, $$ a = a(n,k,r_1,...r_k) = 3^n - 2^{n+k}+\sum_{i=1}^{k} \ 3^{n-r_i}2^{r_i+i-1} $$

\begin{proof} $ $\newline
The lemma will be proven with induction on k.
First of all, it's almost trivial to check that the form of the number will always be $\frac{3^n x_0+a}{2^{n+k}}$, as you multiply by 3, as many times as the horizontal moves you make ($n$) and you divide by 2, as many time as the total moves you make ($n+k$). The whole proof falls into proving the form of the number $a$ that depends on the order, the moves are taken. 
\newline
\newline
We will prove the formula for $a$ in the case $k=0$. 
When $k=0$ we deduce that no $r_i$ exists. So the form of the number will be $$ f(n) = \frac{3^nx_0+a(n)}{2^n}$$ For $n=0$, we have $f(0)=x_0 \Rightarrow a(0)=0$. \\
At the position $(n,k=0)$ the number is of the form given above. Making a horizontal move, at the position $(n+1,0)$ we get $$f(n+1)= \frac{3\frac{3^n x_0+a(n)}{2^n}+1}{2} \Leftrightarrow f(n+1) = \frac{3^{n+1}x_0+3a(n)+2^n}{2^{n+1}} $$ 
From the above equation, we conclude that: $$a(n+1)=3a(n)+2^n , \forall n \in \mathbb{N} $$
The solution of the non-homogeneous linear recurrence is of the form $a_n= a_{HG} + a_{SC}$ where $a_{HG}$ is the solution of the homogeneous recurrence and $a_{SC}$ a special solution. \\
The characteristic equation for the homogeneous recurrence is $r-3=0 \Rightarrow r=3 \Rightarrow a_{HG} =c3^n$, for some constant $c\in R$.\\
We assume a special solution of the form $a_{SC}= b2^n+d$. Substituting this into the equation, gives us $b=-1$ and $d=0$. Using the initial condition $a(0)=0 \Rightarrow c=1$ we receive the general solution: $$ a(n) = 3^n - 2^n, \forall n \in N $$
We notice, that for $k=0$: $$a(n,0)= 3^n-2^{n+0}+\sum_{i=1}^{0} \ 3^{n-r_i}2^{r_i+i-1} = 3^n-2^n $$ and the first case is proven.
\newline
\newline
We assume the formula is correct for k and we prove it for k+1. Let's suppose that the given number $x_0$ has followed the path $(n,k,r_1,...,r_k)$ and has to take the $(k+1)$-th vertical move at $n=r_{k+1}.$ The form of the number before the vertical move will be:

$$ f(x_0,n,k,r_1,...,r_k) = \frac{3^n x_0+a(n,k,r_1,...,r_k)}{2^{n+k}} $$
from the hypothesis. After the vertical move, we get the number: $$\frac{3^n x_0+a(n,k,r_1,...,r_k)}{2^{n+k+1}} $$
But this number must equal to:
$$ f(x_0,n,k+1,r_1,...,r_{k+1})  = \frac{3^nx_0+a(n,k+1,r_1,...,r_{k+1})}{2^{n+k+1}}$$
From the latter equality we extract the relation:
$$ a(n,k+1,r_1,...,r_k,r_{k+1}) = a(n,k,r_1,...,r_k) $$ and by substituting $n=r_{k+1}$ we get:
$$ a(r_{k+1},k+1,r_1,...,r_k,r_{k+1}) = a(r_{k+1},k,r_1,...,r_k) $$ which can be considered as the initial condition.
An application of a horizontal move will give us, as in the first case: $$ \frac{3^{n+1}x_0+3a(n,k+1,r_1,...,r_{k+1})+2^{n+k+1}}{2^{n+k+2}}$$ by which we deduce that $$a(n+1,k+1,r_1,...,r_{k+1})= 3a(n,k+1,r_1,...,r_{k+1}) +2^{n+k+1}$$
It's not hard to prove that the general solution of this non-homogeneous linear recurrence (with respect to $n$) is:
$$ a(n,k+1,r_1,...,r_{k+1}) = c3^n-2^{n+k+1} $$ for some constant c. \\
To calculate c, we put $n=r_{k+1}$ and so we have:
$$ a(r_{k+1},k+1,r_1,...,r_{k+1})= c3^{r_{k+1}}-2^{r_{k+1}+k+1}$$ and by the initial condition:
$$c3^{r_{k+1}}-2^{r_{k+1}+k+1}=a(r_{k+1},k,r_1,...,r_k)= 3^{r_{k+1}}-2^{r_{k+1}+k}+ \sum_{i=1}^{k} \ 3^{r_{k+1}-r_i}2^{r_i+i-1} $$
From the above equation we solve for c to get $$c=1+2^{r_{k+1}+k}3^{-r_{k+1}}+\sum_{i=1}^{k} \ 3^{-r_i}2^{r_i+i-1}$$ and by substituting into the general solution we have:
$$a(n,k+1,r_1,...,r_{k+1}) = (1+2^{r_{k+1}+k}3^{-r_{k+1}}+\sum_{i=1}^{k} \ 3^{-r_i}2^{r_i+i-1})  3^{n}-2^{n+k} =$$
$$ = 3^{n}+3^{n-r_{k+1}}2^{r_{k+1}+k}-2^{n+k+1}+ \sum_{i=1}^{k}\ 3^{n-r_i}2^{r_i+i-1} = 3^n-2^{n+k+1}+3^{n-r_{k+1}} 2^{r_{k+1}+k+1-1}$$
$$ + \sum_{i=1}^{k} \ 3^{n-r_i}2^{r_i+i-1} \Rightarrow a(n,k+1,r_1,...,r_{k+1}) =  3^n-2^{n+k+1} + \sum_{i=1}^{k+1}\ 3^{n-r_i}2^{r_i+i-1} $$  
\end{proof}
$ $\newline
\textbf{Example 2.2.} Let's examine for instance, the path $(n=4,k=3,r_1=0,r_2=2,r_3=4)$. Substituting into the $a$ function gives us: $$a(4,3,0,2,4) = 3^4-2^7+ \sum_{i=1}^{k} \ 3^{4-r_i}2^{r_i+i-1} = -47+3^42^0+3^22^3+3^02^6 = 170 $$ Therefore: $$ x_0 \longmapsto \frac{81x_0+170}{128} $$ If we choose  $x_0 = 118$ we can see that $f(118,4,3,0,2,4) = \frac{81\cdot118+170}{128} = 76$.
Indeed: $$ 118 \Rightarrow 59 \Rightarrow 89 \Rightarrow 134 \Rightarrow 67 \Rightarrow 101 \Rightarrow 152 \Rightarrow 76. $$
\newline
Each path starts from $(n=0,k=0)$. Whether $x_0$ is even or odd the next point on the path will be $(n=0,k=1)$ or $(n=1,k=0)$ respectively. Hence, each path is described by an infinite set of points $(n_i,k_i)$ on a grid. For example, the path $(x_0,n=2,k=2,r_1=0,r_2=2)$ is given by the path sequence:
 $$ (0,0) \xrightarrow[\text{}]{\text{$r_1=0$}} (0,1) \xrightarrow[]{} (1,1) \xrightarrow[]{} (2,1) \xrightarrow[\text{}]{\text{$r_2=2$}} (2,2) $$
The next lemma about $f$ is crucial in the discovery of
the formula for the characteristic function.
\newline
\newline
\textbf{Lemma 2.3}
Consider the path $(n,k,r_1,...,r_k)$ that $x_0$ follows and its corresponding path sequence. Then, for every value of $n_0$ in the path sequence, there exists exactly one $k_0$ for which $f(x_0,n_0,k_0,r_1,...r_{k_0})$ is odd. Equally:
$$\forall n, \exists! k : f(x_0,n,k,r_1,...,r_k)=2\lambda_{n}+1$$
\begin{proof} $ $\newline
First notice that more than one point in a sequence may have the same $n$. For instance, in a vertical move we go from $(n,k)$ to $(n,k+1)$. We fix a single $n_0$ and iterate through all points that have $n=n_0$. If for the first point $(n_0,k_0), f(x_0,n_0,k_0,r_1,...,r_{k_0})$ is odd, then the next point must be $(n_0+1,k_0)$ and thus the only $k$ satisfying the lemma for $n_0$ is $k_0$. If $f(x_0,n_0,k_0,r_1,...,r_{k_0})$ is even, then the next point must be $(n_0,k_1=k_0+1)$ and so, $k_0$ does not satisfy the lemma. This process is repeated until $f(x_0,n_0,k_i,r_1,...,r_{k_i})$ is odd, for some $k_i=k_0+i$ which will be the unique number satisfying the lemma for $n_0$.
\end{proof}
$ $\newline
The lemma equation can be equally written as
$$ \frac{3^nx_0+a}{2^{n+k}}=2\lambda_n+1$$
The above relation is true for any n and for these k such that $f(n,k,r_1,...,r_k)$ is odd. Hence, if we solve this equation for $x_0$ in terms of $n,k,r_1,...,r_k$ we have found the characteristic equation. It is
$$ 3^nx_0 + 3^n-2^{n+k} +  \sum_{i=1}^{k}\ 3^{n-r_i}2^{r_i+i-1} = 2^{n+k+1}\lambda_n+2^{n+k} $$
or
$$ -3^n(x_0+1) + 2^{n+k+1}(\lambda_n+1) = \sum_{i=1}^{k}\ 3^{n-r_i}2^{r_i+i-1}$$
\textbf{Theorem 2.3. The Characteristic Function} 
\newline
The solution to the Diophantine equation
$$ -3^n(x_0+1) + 2^{n+k+1}(\lambda_n+1) = \sum_{i=1}^{k}\ 3^{n-r_i}2^{r_i+i-1}$$
is given by:
{\Large $$\begin{cases}  x_0 = 2^{n+k+1}m -1-\sum_{i=1}^{k}\ \frac{c_i2^{n+k+1}+2^{r_i+i-1}}{3^{r_i}} \\ \lambda_n=3^nm-1-\sum_{i=1}^{k}\ c_i3^{n-r_i} \end{cases} $$ }
where $c_i$ is a family of natural constants and the equation holds $\forall m \in \mathbb{Z} $
\begin{proof} $ $\newline
The Diophantine equation $ax+by=c$, with $a,b,c \in \mathbb{Z}$ has a solution if and only if $gcd(a,b)|d$. Moreover, given a special solution $(x^*,y^*)$ the general solution of the equation is of the form: 
$$ (x,y)=\Big( \frac{b}{gcd(a,b)}m + x^*,-\frac{a}{gcd(a,b)}m+y^*\Big), \forall m \in Z $$ 
This above result is known from the theory of Diophantine equations. The to-be-solved equation meets the criteria of a Diophantine equation because $-3^n,2^{n+k+1}$ and $\sum_{i=1}^{k}\ 3^{n-r_i}2^{r_i+i-1}$
are integer numbers ($r_i\leq n).$ \\
It can easily be proven that the equation has infinitely many solutions. Indeed, $gcd(a,b)=gcd(-3^n,2^{n+k+1})=1$ which divides any integer number. So $1|\sum_{i=1}^{k}\ 3^{n-r_i}2^{r_i+i-1}$ and thus:
$$ x_0+1 = \frac{b}{gcd(a,b)}m+x^* = 2^{n+k+1}m+x^* $$ and 
$$ \lambda_n+1= -\frac{a}{gcd(a,b)}m+\lambda^* = 3^nm+\lambda^* $$
So, if the special solution: $$(x^*,\lambda^*)= \Big( -\sum_{i=1}^{k}\ \frac{c_i2^{n+k+1}+2^{r_i+i-1}}{3^{r_i}}, -\sum_{i=1}^{k}\ c_i3^{n-r_i}\Big) $$
verifies the given equation and moreover $(x^*,\lambda^*) \in {\mathbb{Z}}^2$ the lemma is proved.
\newline
In fact, it's just a matter of operations to show that the special solution satisfies the equation. $$ -3^n( -\sum_{i=1}^{k}\ \frac{c_i2^{n+k+1}+2^{r_i+i-1}}{3^{r_i}}) + 2^{n+k+1}(-\sum_{i=1}^{k}\ c_i3^{n-r_i}) = $$ 
$$ \sum_{i=1}^{k}\ c_i3^{n-r_i}2^{n+k+1} + \sum_{i=1}^{k}\ 3^{n-r_i}2^{r_i+i-1}- \sum_{i=1}^{k}\ c_i3^{n-r_i}2^{n+k+1} = \sum_{i=1}^{k}\ 3^{n-r_i}2^{r_i+i-1} $$
\\
It's obvious that for $c_i \in \mathbb{Z} \Rightarrow \lambda^* \in \mathbb{Z}$. For $x^*$ we need to show that $$\forall i \in \{1,...,k\}, \exists c_i \in \mathbb{Z} : \frac{c_i2^{n+k+1}+2^{r_i+i-1}}{3^{r_i}} \in \mathbb{Z}.$$
Equivalently, we can write:
$$ c_i2^{n+k+1}+2^{r_i+i-1} = p_i3^{r_i} \Leftrightarrow -2^{n+k+1}c_i + 3^{r_i}p_i = 2^{r_i+i-1},\text{for } p_i \in \mathbb{Z}. $$
The existence of $c_i$ is secured by the fact that $gcd(-2^{n+k+1},3^{r_i}) = 1 | 2^{r_i+i-1}, \forall i$, since the above is a Diophantine equation. 
\end{proof}
$ $\newline
Through a couple of examples, the characteristic function will become clearer of how it works.
\newline
\newline
\textbf{Example 2.4.} Consider the following path $(n=4,k=2,r_1=3,r_2=4)$. We are looking forward to finding the \textbf{smallest} number that follows this path. By substituting into the characteristic function we get:
$$ x(4,2,3,2) = 2^7m-1-\sum_{i=1}^{2}\ {\frac{2^7+2^{r_i+i-1}}{3^{r_i}}} \Rightarrow$$
$$ x = 128m-1-\frac{128c_1+2^{3+1-1}}{3^3}-\frac{128c_2+2^{4+2-1}}{3^4}=128m-1-\frac{128c_1+8}{27}-\frac{128c_2+32}{81}. $$
\\
There exist $c_1,c_2$ such that the two fractions are integer numbers. By trying out numbers between $0$ and $27$ we can see that a (special) solution to this hidden form of a Diophantine equation is $c_1=5$ and likewise $c_2=20$.
\newline
So, the general form of $c_1,c_2$ will be:
$$ c_1 = 27m_1+5, \text{ } c_2 = 81m_2+20 ,\text{ } \forall m_1,m_2 \in \mathbb{Z}. $$
Substituting back to $x$ we get:
$$ x=128m-1-128m_1-128m_2-\frac{128\cdot5+8}{27}- \frac{128\cdot20+32}{81} \rightarrow $$  
$$ x=128m_*-1-24-32 \rightarrow  x = 128m-57   \xrightarrow[\text{}]{\text{$m=1$}}x_{min} = 71. $$
We should check to see if the number 71 follows this path. Indeed: 
$$ 71 \xrightarrow[\text{}]{\text{$n=1$}} 107 \xrightarrow[\text{}]{\text{$n=2$}} 161 \xrightarrow[\text{}]{\text{$n=3$}} 242 \xrightarrow[\text{$r_1=3$}]{\text{$k=1$}}121
\xrightarrow[\text{}]{\text{$n=4$}} 182
\xrightarrow[\text{$r_2=4$}]{\text{$k=2$}} 91. $$
\newline
\textbf{Example 2.5.}
Let's now consider, the simpler path $(n=1,k=1,r_1=0)$. Without any work, we can see that the number x will be even, in order to run a vertical move instantly. By substituting into the function we get:
$$ x(1,1,0) = 2^3m-1-\frac{2^3c_1+2^{0+1-1}}{3^0} = 8m-1 - (8c+1) = 8m-2. $$
So, the smallest number that follows this path is, for $m=1$, $x=6$. We can easily check that:
$$ 6 \Rightarrow 3\Rightarrow 5.$$
However, one can see that $x=2$ also follows the path $(n=1,k=1,r_1=0)$ ($2\Rightarrow1 \Rightarrow2$) and is smaller than 6. The characteristic function is not incorrect. Examining the way the function was built, it can be seen that given a path, the function returns the smallest number that follows this path but results in an odd number. On our specific example, if we need to find the answer $x=2$, we should give as input the path $(n=1,k=2,r_1=0,r_2=1)$ and that's only because the next move will be horizontal. We note here, that the meaning of the function is to find the number that follows a specific path and not to find what path should we give to the function in order to get a specific number. 
\newline
\newline
\textbf{Example 2.6.} Consider the path $(n,k=n,r_1=1,...,r_k=k)$. This complicated written, path isn't other than the "one horizontal move, one vertical move, one horizontal move, one vertical move, ...". It's not hard to guess which number is the smallest one to follow this path but let's try to prove it. By substituting we get:
$$ x = 2^{2n+1}m-1-\sum_{i=1}^{n}\ \frac{c_i2^{2n+1}+2^{2i-1}}{3^i} = 2^{2n+1}m-1-2^{2n+1}\sum_{i=0}^{n}\ \frac{c_i}{3^i} - \sum_{i=1}^{n}\ \frac{2^{2i-1}}{3^{i}}. $$
\newline
Although, we separated the fraction into two terms, we must not forget that c must be chosen such that $x\in \mathbb{Z}$. 
\newline
The second sum is known and can be easily computed as follows:
$$ \sum_{i=1}^{n}\ \frac{2^{2i-1}}{3^{i}} = \frac{1}{2} \sum_{i=1}^{n}\ (\frac{4}{3})^i = 2[(\frac{4}{3})^n-1] = \frac{2^{2n+1}}{3^n}-2 $$
\newline
For the first sum, as we have no clue on what $c_i$ could be, we assume that $c_i=c \text{ },\forall i$. If we cannot find $c$ such that $x\in \mathbb{Z}$ then our assumption was wrong. So the first sum will be equal to:
$$ \sum_{i=0}^{n}\ \frac{c_i}{3^i} = c\sum_{i=0}^{n}\ (\frac{1}{3})^i = \frac{c}{2} (1-\frac{1}{3^n}) = \frac{c}{2} - \frac{1}{2} \frac{c}{3^n}. $$
\newline
So, we substitute the result back on our first equation. We have:
$$ x= 2^{2n+1}m-1-c2^{2n}+c\frac{2^{2n}}{3^n}-\frac{2^{2n+1}}{3^n}+2. $$
\newline
We observe that for $c=2$ the two fractions cancel out, leaving only integer terms in x. Now we can, finally, calculate the number x:
$$ x = 2^{2n+1}m-1-2^{2n+1}+2 \xrightarrow[\text{}]{\text{$m=1$}} x=1. $$
\newline
Note that $(n,n,r_i=i)$ is the \textbf{characteristic path} for the number 1. Generally, it is significantly hard to prove a characteristic path $(n,k=f(n),r_i=g(i))$ for a number. For instance, in our last example, if $c_i$ where not all the same, the proof would become a lot harder.
\newline
\newline
Having found and understood the characteristic function,we head into the last part of the proof of the Collatz Conjecture.
\newpage
\section{The Proof}
Recall the function $C[x]$ from Section 1 and let $x_0$ be the smallest number that never reaches 1. In other words $C[x_0] = 0$ and 
$$\text{} C[x]=1,  \forall x\in \{1,2,...,x_0-1\} $$
Given a sequence $x_0,x_1,...,x_n,...,$ if $C[x_i]=1$  for any $i>0$ then $C[x_0] = 1$ as well.
\newline
\newline
Supposing that the smallest number that never reaches 1, namely $x_0$, is even, then $x_0 \rightarrow x_1=\frac{x_0}{2}$. But $\frac{x_0}{2}<x_0$ which means that $C[\frac{x_0}{2}]=1$, concluding that $C[x_0]=1$, which is a contradiction to our hypothesis. So $x_0$ cannot be an even number.
\newline
\newline
To extend the previous example, suppose that $x_0$ follows the path $(n=1,k=2,r_1=1,r_2=1)$. We can compute the outcome number at the end of the path to be:
$$ x_3 = f(x_0,n=1,k=2,r_1=1,r_2=1) = \frac{3^1x_0+a(1,2,1,1)}{2^{1+2}} = \frac{3x_0+1}{8} <x_0. $$
Likewise, that means that $C[x_0]=1$, which is a contradiction. So, $x_0$ can't follow the path $(n=1,k=2,r_1=1,r_2=1)$. The next lemma generalizes the former examples.
\newline
\newline
\textbf{Lemma 3.1} Let $x_0$ be a positive integer and $(n,k,r_1,...,r_k)$ a path that $x_0$ follows. Then $f(x_0,n,k,r_1,...,r_k)<x_0$ if and only if $2^{n+k}>3^n$.
\begin{proof} $ $\newline
The -to be proven- inequality is written as:
$$ \frac{3^nx_0+a(n,k,r_1,...,r_k)}{2^{n+k}} < x_0 \Leftrightarrow (2^{n+k}-3^n)x_0 > a(n,k,r_1,...,r_k). $$
If $2^{n+k}\leq 3^n$ then the right hand side is a negative number, unable to be greater than a positive number, namely the right hand side. So, in order the inequality to hold $2^{n+k}>3^n$.
\newline
We will prove now that $$x_0> \frac{a(n,k,r_1,...,r_k)}{2^{n+k}-3^{n}}$$ while $2^{n+k}>3^n$. Suppose $x_0$ is the smallest number that follows the path $(n,k,r_1,...,r_k)$. If we prove the inequality for the smallest number, then any bigger number that follows the same path will satisfy the inequality since
$$ x^*_0 > x_0 > \frac{a(n,k,r_1,...,r_k)}{2^{n+k}-3^{n}}$$ 
So $x_0$ is given by the characteristic function $F(n,k,r_1,...,r_k)$. For that to be the case though, $f(x_0,n,k,r_1,...,r_k)$ must be odd since $F$ returns the smallest number that follows a path that ends to an odd number (Example 2.5). Suppose that $f$ is even. Then we can extend $x_0$'s path to apply a vertical move and the path will become $(n,k+1,r_1,...,r_k,r_{k+1}=n)$ and $f(x_0,n,k+1,r_1,...,r_{k+1}) = \frac{1}{2}f(x_0,n,k,r_1,...,r_k).$ But $2^{n+k+1}>2^{n+k}>3^n$ and so the inequality's proof can be reduced to
 $$ \frac{3^n x_0+a(n,k,r_1,...,r_{k+1})}{2^{n+k+1}}=\frac{f(x_0,n,k,r_1,...,r_k)}{2} < x_0 $$ where $k' = k+1$. If $f(x_0,n,k+1,r_1,...,r_{k+1})$ is even, the process can be repeated until $f(n,k+l,r_1,...,r_{k+l})$ is odd or
 $$\frac{f(x_0,n,k,r_1,...,r_k)}{2^l} < x_0 $$
 for some positive integer $l$.
 Hence, we can assume that $f(x_0,n,k,r_1,...,r_k)$ is odd and that $x_0$ is given by the characteristic function $F(n,k,r_1,...,r_k)$. The inequality then, can be equally written:
 $$ 2^{n+k+1}m-1-\sum_{i=1}^{k}\ \frac{c_i2^{n+k+1}+2^{r_i+i-1}} {3^{r_i}}>\frac{3^n-2^{n+k}+\sum_{i=1}^{i=k}\ 3^{n-r_i}2^{r_i+i-1}}{2^{n+k}-3^n}\Leftrightarrow$$
 $$ \Leftrightarrow 2^{2n+2k+1}m-2^{2n+2k+1}\sum_{i=1}^{k} \ \frac{c_i}{3^{r_i}}-2^{n+k}\sum_{i=1}^{k}\ \frac{2^{r_i+i-1}} {3^{r_i}}-3^n2^{n+k+1}m+2^{n+k+1}\sum_{i=1}^{k}\ c_i3^{n-r_i}>0 \Leftrightarrow $$
$$ \Leftrightarrow 2^{n+k+1}m-2^{n+k+1}\sum_{i=1}^{k} \ \frac{c_i}{3^{r_i}}-2\cdot 3^nm+2\sum_{i=1}^{k}\ c_i3^{n-r_i}> \sum_{i=1}^{k}\ \frac{2^{r_i+i-1}} {3^{r_i}} \Leftrightarrow $$ $$\Leftrightarrow
2^{n+k+1}(m-\sum_{i=1}^{k} \ \frac{c_i}{3^{r_i}})-2\cdot3^n(m -\sum_{i=1}^{k} \ \frac{c_i}{3^{r_i}}) > \sum_{i=1}^{k}\ \frac{2^{r_i+i-1}} {3^{r_i}} \Leftrightarrow$$
$$ \Leftrightarrow 2(2^{n+k}-3^n)(m-\sum_{i=1}^{k} \ \frac{c_i}{3^{r_i}})>\sum_{i=1}^{k}\ \frac{2^{r_i+i-1}} {3^{r_i}} $$
In order to prove the last inequality we make two observations. 
\begin{itemize}
    \item We choose $m$ and $c_i$ so that $m-\sum_{i=1}^{k} \ \frac{c_i}{3^{r_i}} > 1$. In that way we can say that: $$ 2(2^{n+k}-3^n)(m-\sum_{i=1}^{k} \ \frac{c_i}{3^{r_i}}) > 2(2^{n+k}-3^n).$$
    \item It is $r_i\geq 0.$ Therefore, 
    $$(\frac{2}{3})^{r_i} \leq 1 \Rightarrow \sum_{i=1}^{k}\ \frac{2^{r_i+i-1}} {3^{r_i}} \leq \sum_{i=1}^{k}\ 2^{i-1}= 2^k-1. $$
\end{itemize}
So, we get: 
$$ 2(2^{n+k}-3^n)(m-\sum_{i=1}^{k}\ \frac{c_i}{3^{r_i}}) > 
2(2^{n+k}-3^n) \stackrel{?}{>} 2^k-1 \geq  \sum_{i=1}^{k}\ \frac{2^{r_i+i-1}} {3^{r_i}}. $$
The -to be proven- inequality is written equally:
$$ 2^{n+1}2^k-2\cdot 3^n > 2^k-1 \Leftrightarrow 2^k(2^{n+1}- 1) > 2\cdot 3^n-1 \Leftrightarrow 2^k> \frac{2\cdot 3^n-1}{2 \cdot2^n -1}.$$
We have supposed though, that $2^{n+k}>3^n \Rightarrow 2^k > (\frac{3}{2})^n$. So, if we prove that:
$$ 2^k > \frac{3^n}{2^n} \Longrightarrow 2^k >\frac{3^n- \frac{1}{2}}{2^n-\frac{1}{2}}$$
then, the initial inequality is proven. 
\newline
\newline
We can equally write the last inequality as:
$$ \frac{2^{k+n}}{3^n} > \frac{1-\frac{1}{2\cdot3^n}}{1-\frac{1} {2\cdot2^n}} \Leftrightarrow (k+n)log2-nlog3 > log(1-\frac{1}{2\cdot3^n})-log(1-\frac{1} {2\cdot2^n}). $$
G.Rhin has proven that: 
$$ (k+n)log2-nlog3 > \frac{1}{457\cdot n^{13.3}} $$
for every positive integer $n,k$. \textbf{[1],[2]}
We can easily now verify, through Wolfram for instance,  that for every $n\geq 96$:
$$ \frac{1}{457\cdot n^{13.3}} >log(1-\frac{1}{2\cdot 3^n})-log(1-\frac{1}{2\cdot 2^n}) $$
and so the initial inequality is proven. For the remaining values of $n$ we can check through a computer that the initial inequality holds.
\end{proof}
$ $\newline
It has been proven that a number $x_0$ that follows the path $(n,k,r_1,...,r_k)$ with $2^{n+k}>3^n$ will satisfy the inequality:
$$ x_{n+k} = f(x_0,n,k,r_1,...,r_k) = \frac{3^nx_0+a(n,k,r_1,...,r_k)}{2^{n+k}}<x_0$$
Thus, if $x_0$ is the smallest number that $C[x_0]=0$ then $C[x_{n+k}] = 1$ since $x_{n+k}<x_0$ implying that $C[x_0]=1$, a contradiction. So $x_0$ cannot follow any path with $n$ and $k$ satisfying $2^{n+k}>3^n$.
\newline
\newline
The inequality can be solved for n:
$$ 2^{n+k}>3^n \Leftrightarrow log_22^{n+k}>log_23^n \Leftrightarrow n+k > nlog_23 \Leftrightarrow n< \frac{k}{log_23-1}. $$
We can use the last relation to deduce a restriction for $r_i$'s. Specifically $r_i\leq n$ and replacing $k$ with $i$ yields:
$$r_i < \frac{i}{log_23-1}.$$
Some values for $r_i$'s are $r_1<2, r_2<4, r_3<6$ and so on. If $x_0$ follows a path with at least one $r_i$ satisfying the above inequality, then $x_0$ will eventually reach 1. Summing up, in order for $x_0$ to never reach 1, namely $C[x_0]=0$, it must follow a path $(n,k,r_1,...,r_k)$ where:
 $$\begin{cases}  r_i>\frac{i}{log_23-1}, \forall i\in \mathbb{N} \\ k<n(log_23-1) \end{cases} $$
Before heading to the last theorem that completes the proof of the Collatz Conjecture,we introduce a final example.
\newline
\newline
\textbf{Example 3.2} Suppose we want to find the smallest number that makes only horizontal moves. The path for this will be $(n,0)$ with $n\rightarrow \infty.$ From the characteristic function we can easily find  $x_0$ to be:
$$x_0=2^{n+k+1}m-1-\sum_{i=1}^{k}\ \frac{c_i2^{n+k+1}+2^{r_i+i-1}}{3^{r_i}}= 2^{n+1}m-1 \Rightarrow x_{min} = 2^{n+1}-1.$$
So, to find the smallest number that does 4 horizontal moves, for example, we substitute $n=4$ in the above equation and get $x=2^{5}-1=31$. Indeed:
$$ 31 \Rightarrow 47 \Rightarrow 71 \Rightarrow 107 \Rightarrow 161. $$
The smallest number that does infinite many horizontal moves, without doing a single vertical one, and consequently never reaching 1,  will be:
$$ x_0 = \lim_{n\to\infty}x_{min} = \lim_{n\to\infty} (2^{n+1}-1) = +\infty. $$
We conclude that there is no natural number to run horizontal moves forever. 
\newline
\newline
We move straight to the last lemma.
\newline
\newline
\textbf{Lemma 3.3} There exists no natural number $x_0$ that follows the path $(n,k,r_1,...,r_k)$ with
$$\begin{cases}  r_i>\frac{i}{log_23-1}, \forall i\in \mathbb{N} \\ k<n(log_23-1) \end{cases} $$
\begin{proof} $ $\newline
First, notice that $$ (\frac{2}{3})^{\frac{i}{log_23-1}} = 2^{-i} $$
Indeed, 
$$  (\frac{2}{3})^{\frac{i}{log_23-1}} = 2^{-i} \Leftrightarrow 2^i2^{\frac{i}{log_23-1}}= 3^{\frac{i}{log_23-1}} \Leftrightarrow 2^{\frac{ilog_23} {log_23-1}} = 3^{\frac{i}{log_23-1}}\Leftrightarrow $$
$$ \Leftrightarrow 2^{ilog_23}=3^i \Leftrightarrow (2^{log2_3})^{i} = 3^i \Leftrightarrow 3^i = 3^i$$
which is true.
Now, $x_0$ will be given by the characteristic function
$$x_0=2^{n+k+1}m-1-\sum_{i=1}^{k}\ \frac{c_i2^{n+k+1}+2^{r_i+i-1}}{3^{r_i}} = 2^{n+k+1}(m-\sum_{i=1}^{k}\ \frac{c_i}{3^{r_i}})-1- \sum_{i=1}^{k}\ (\frac{2}{3})^{r_i}2^{i-1}.$$
Now, we can notice that for $r_i > \frac{i}{log_23-1}$:
$$ \sum_{i=1}^{k}\ (\frac{2}{3})^{r_i}2^{i-1} < \sum_{i=1}^{k}\ (\frac{2}{3})^{\frac{i}{log_23-1}}2^{i-1} = \sum_{i=1}^{k}\ 2^{-i}2^{i-1} = \sum_{i=1}^{k}\ \frac{1}{2} = \frac{k}{2}. $$
And by using the fact that $k<n(log_23-1)$ we get:
$$\sum_{i=1}^{k}\ (\frac{2}{3})^{r_i}2^{i-1} < \frac{n(log_23-1)}{2}.$$
Hence, the smallest number that follows such path continuously is given by the limit
$$ x_0=\lim_{n\to\infty}\Big[2^{n+k+1}(m-\sum_{i=1}^{k}\ \frac{c_i}{3^{r_i}})-1- \sum_{i=1}^{k}\ (\frac{2}{3})^{r_i}2^{i-1}\Big] $$
Lastly, $m$ and $c_i$ are chosen so $x_0$ is the smallest number that follows the path and obviously
$$m-\sum_{i=1}^{k}\ \frac{c_i}{3^{r_i}}>0$$
a positive constant number. Finally,
$$ x_0 = \lim_{n\to\infty}\Big[2^{n+k+1}(m-\sum_{i=1}^{k}\ \frac{c_i}{3^{r_i}})-1- \sum_{i=1}^{k}\ (\frac{2}{3})^{r_i}2^{i-1}\Big] > \lim_{n\to\infty} \Big(2^{n+k+1}m_*-1-n\frac{log_23-1}{2}\Big)= +\infty $$
That means:
$$ x_0 > +\infty \Rightarrow x_0 = +\infty$$
and so there exists no natural number that follows a path $(n\to\infty, k\to\infty,r_1,...,r_k,...)$ that infinitely satisfies
$$\begin{cases}  r_i>\frac{i}{log_23-1}, \forall i\in \mathbb{N} \\ k<n(log_23-1) \end{cases} $$
\newline
\end{proof}
$ $\newline
\textbf{Theorem 3.4} There exists no natural number $x_0$ for which $C[x_0]=0$.
\begin{proof} $ $\newline
By lemma 3.3 it has been proven that a natural number $x_0$ cannot follow a path with $r_i>\frac{i}{log_23-1}, \forall i\in \mathbb{N}$ and $2^{n+k}>3^n$ infinitely. Thus, it must follow a path in which at least one move satisfies $2^{n+k}<3^n$. By lemma 3.1, though, in this move, number $x_{n+k}=f(x_0,n,k,r_1,...,r_k)<x_0$ concluding that $C[x_{n+k}]=1$. Consequently $C[x_0]=1$ for any natural number. 
\newline
Notice, that in the path $x_0,x_1,...,x_{n+k}$, there is no restrain for the numbers to be different. In Example 2.6, the path $(n,k,r_1=1,...,r_k=k)$ is $1,2,1,2,...$. If a cycle $x_0,...,x_{n+k},x_0$ existed such that in the infinite path,  $r_i>\frac{i}{log_23-1}$ was always true, that would contradict lemma 3.3 and by that no cycle can exist without reaching 1.  
\end{proof}
$ $\newline

\vspace*{\fill}
\noindent \textbf{Agelos Kratimenos} \newline
National Technological University of Athens (NTUA) \newline
Athens, Greece \newline
E-mail: ageloskrat@yahoo.gr
\end{document}